\documentclass{article}
\usepackage{amssymb,amsmath,graphicx}
\begin{document}

\title{Perfect Algebraic Coarsening}
\author{Jonathan E. Moussa\thanks{Department of Physics, University of California at Berkeley, Berkeley, California 94720. Materials Science Division, Lawrence Berkeley National Laboratory, Berkeley, California 94720. \textit{E-mail address:} jmoussa@civet.berkeley.edu. }} 

\maketitle

\begin{abstract}
Presented in this paper is a new sparse linear solver methodology
 motivated by multigrid principles and
 based around general local transformations
 that diagonalize a matrix while maintaining its sparsity.
These transformations are approximate,
 but the error they introduce can be systematically reduced.
The cost of each transformation is independent of matrix size
 but dependent on the desired accuracy
 and a spatial error decay rate governed by local properties of the matrix.
We test our method by applying a single transformation
 to the 2D Helmholtz equation at various frequencies,
 which illustrates the success of this approach.
\end{abstract}




\section{Introduction}

Physicists develop mathematical models from physical considerations,
 but the process of solving a model isn't always related to its physics.
Intermediate steps of a long calculation may not have physical meaning
 nor grant physical insight.
An important example of this is sparse linear system solving,
 which is used to solve discretized approximations of physical problems
 described by time-independent partial differential equations.
The most general-purpose sparse linear solvers involve either direct factorization \cite{Golub},
 where intermediate steps contain partially factored matrices,
 or iterative solvers \cite{templates},
 where intermediate steps contain approximate solutions of increasing accuracy.
For some well-understood problems,
 most notably the Poisson equation,
 known physical properties can be incorporated into a linear solver,
 either via multigrid methods \cite{multigrid} or multilevel preconditioning \cite{multilevel},
 leading to algorithms that are both more physical and of optimal complexity.
These methods operate by approximately transforming away local details of the physical system,
 leaving successively smaller
 but continually sparse ``coarsened" matrix equations
 that each represent the physical system on a different length scale.
The multigrid framework \cite{generalAMG} and coarsening procedures \cite{AMGcoarsening}
 have been generalized into a more algebraic formalism,
 but their success is still tied to certain spectral properties of the underlying physical system.
In this paper we construct a more general sparse matrix solver
 based on a high accuracy limit of algebraic multigrid
 that does not require as input a physical intuition for the problem
 and is not restricted by spectral properties.

The test system we use in this paper
 is the 2D Helmholtz equation
 discretized on a uniform grid using finite differences
 and defined by the five point matrix stencil
\begin{displaymath}\left[ \begin{array}{ccc}  & 1 &  \\ 1 & (\lambda-4) & 1 \\  & 1 &  \end{array} \right], \end{displaymath}
 where $\lambda$ is proportional to the frequency squared.
The finite difference approximation loses accuracy as $\lambda$ gets larger
 and breaks down entirely for $\lambda > 4$ due to inadequate sampling of oscillations,
 but we are interested in the matrix problem
 and not necessarily its accuracy in reproducing the continuum problem.
This is perhaps the simplest example
 of a problem for which optimal linear solvers exist
 but as yet require some analytic knowledge of the solution to construct.
For rectangular domains,
 the eigenfunctions are composed of sinusoidal oscillations,
 and thus fast fourier transforms \cite{FFT} can diagonalize the matrix.
Given the exact analytic inverse of the Helmholtz equation,
 one can hierarchically compress and apply it in an optimal manner using the fast multipole method \cite{FMM}.
In a more algebraic manner,
 using just the knowledge that solutions have a characteristic frequency of oscillation,
 it is possible to construct an optimal ray-based multigrid scheme \cite{waveAMG}.
Algebraic methods that don't take specific account of the oscillatory nature of solutions
 currently fail to solve the problem in an optimal manner.
For direct factorizations,
 the cost has been proven to be non-optimal for problems on a 2D grid \cite{nesteddissection}.
For the parameter range of oscillatory behavior, $0 < \lambda < 8$,
 preconditioners based on multigrid principles fail to be optimal
 due to a loss of smoothness on coarse grids \cite{helmiter},
 and structured direct methods fail due to the loss of low off-diagonal rank \cite{hackbusch}.

In Section \ref{formtrans}, we describe the form of the linear solver
 as a succession of local transformations
 and some of their properties and governing equations.
In Section \ref{linmethod}, we derive an efficient method for constructing local transformations
 and apply it to our test problem.
In Section \ref{blocksection}, we further generalize the local transformations
 by changing sparsity patterns to improve accuracy.

\section{Form of the Transformation\label{formtrans}}

The basic operation of our linear solver is to start from $A$,
 an $n \times n$ real symmetric sparse matrix at some stage of factorization,
 and apply a real symmetric transformation,
\begin{displaymath} X^T A X = \tilde{A} + E, \end{displaymath}
 that leaves us with an $\tilde{A}$ that has one more diagonalized row/column than $A$
 and a small amount of error, $E$.
Transformations of this form are found in direct LDL factorization \cite{Golub},
 where $X$ is the identity plus a rank one matrix
 and $E$ is just floating point roundoff error.
This form can also be related to multigrid solvers,
 if coarsening and relaxation are combined into a single invertible transformation \cite{generalAMG}
 and if coarsening is performed only on one single small subdomain at a time.
The $X$ will be a dense matrix on this small subdomain and identity outside of it,
 and $E$ on the subdomain will be substantially larger than roundoff errors.

The benefit of multigrid,
 despite the large error,
 is that the sparsity pattern of $\tilde{A}$ can be more controlled
 and the critical filling in of $A$ during factorization
 that prevents LDL factorization from being optimal can be avoided.
The large error $E$ incurred at each factorization step can be negated
 by including a multilevel refinement scheme in the linear solving procedure,
 but the success of refinement is based on
 details of the spectral properties of the problem \cite{generalAMG}.
The current prescription for reducing multigrid coarsening error
 is simply to increase fill-in of the sparsity pattern of $\tilde{A}$ \cite{AMGcoarsening},
 but this relates error to matrix fill
 and reduces our ability to control the sparsity pattern.
Our more general approach is to hold the sparsity pattern of $\tilde{A}$ fixed
 while allowing more freedom in the choice of $X$,
 enough to enable $\|E\|$ to be arbitrarily reduced,
 bounded only by machine precision.
We denote this fixed sparsity, high accuracy limit of algebraic coarsening
 as \textit{perfect algebraic coarsening}
 and denote the $X$ matrices as \textit{local sparsity-preserving transformations}.

Repeated transformations take us
 from our initial matrix $A_0$
 to a final diagonal form $D$,
\begin{displaymath} X_n^T \cdots X_1^T A_0 X_1 \cdots X_n = D + O(\|E\|), \end{displaymath}
 which leads to a compact representation of the inverse of $A_0$,
\begin{displaymath} A_0^{-1} = X_1 \cdots X_n D^{-1} X_n^T \cdots X_1^T + O(\|E\|). \end{displaymath}
If $\|E\|$ can be reduced sufficiently,
 then this multigrid-based factorization can be made as accurate as a direct factorization,
 foregoing the need for the iterative steps of multigrid.
If we can restrict each transformation $X_i$
 to differ from identity only on an $n$-independent sized subdomain of the problem,
 then each of these $n$ transformations can be calculated with an $n$-independent cost,
 and the resulting linear solver will be of optimal $O(n)$ complexity.

The restriction on each $X$ is a ``local" one,
 which in terms of the underlying grid
 means that a transformation that decouples a node
 should only act on neighbors of that node
 up to at most some $m^\textrm{th}$ nearest neighbor.
The restricted transformation takes the form
\begin{equation} \left[ \begin{array}{cc} X_L^T & 0 \\ 0 & I \end{array} \right] \left[ \begin{array}{cc} A_{LL} & A_{LE} \\ A_{LE}^T & A_{EE} \end{array} \right] \left[ \begin{array}{cc} X_L & 0 \\ 0 & I \end{array} \right] = \left[ \begin{array}{cc} \tilde{A}_{LL} & A_{LE} \\ A_{LE}^T & A_{EE} \end{array} \right] + E, \label{localtransform}\end{equation}
 where the subscript `$L$' refers to a local partition
 and `$E$' to the remaining external partition.
In order for Eq. (\ref{localtransform}) to be satisfied with a small error,
 we have to enforce the condition, $A_{LE}^T X_L = A_{LE}^T$,
 either approximately with some least squares approach
 or exactly by finding the null space of $A_{LE}^T$
 or more simply by further partitioning the local region
 into an interior `$I$' and boundary `$B$',
\begin{displaymath} A \rightarrow \left[ \begin{array}{ccc} A_{II} & A_{IB} & 0 \\ A_{IB}^T & A_{BB} & A_{BE} \\ 0 & A_{BE}^T & A_{EE} \end{array} \right], \end{displaymath}
 and further restricting $X_L$ to
\begin{equation} \left[ \begin{array}{cc} X_I^T & 0 \\ X_B^T & I \end{array} \right] \left[ \begin{array}{cc} A_{II} & A_{IB} \\ A_{IB}^T & A_{BB} \end{array} \right] \left[ \begin{array}{cc} X_I & X_B \\ 0 & I \end{array} \right] = \left[ \begin{array}{cc} \tilde{A}_{II} & \tilde{A}_{IB} \\ \tilde{A}_{IB}^T & \tilde{A}_{BB} \end{array} \right] + E. \label{interiortransform}\end{equation}
Our calculation of $X$ and $\tilde{A}$ may now proceed independently of the external partition,
 with some $n$-independent cost dependent only on
 $n_L$, $n_B$, and $n_I$ -- the sizes of the local, boundary, and internal partitions --
 and $n_{\tilde{A}}$, the number of independent nonzero elements in the symmetric $\tilde{A}_{LL}$.

We must next define an error norm to be minimized by our choice of transformation.
A convenient choice of norms when dealing with variable matrices
 is the Frobenius norm, $\|M\|_F \equiv \sqrt{Tr(M^H M)}$.
Minimizing $\|E\|_F$ directly leads to the expression
\begin{equation} \min_{X,\tilde{A}} {\left\| X^T A X - \tilde{A} \right\|_F}, \label{errmin}\end{equation}
 with $\tilde{A}$ restricted to a given sparsity pattern
 and $X$ restricted to the form in Eq. (\ref{interiortransform}).
This error norm is problematic because it is dependent on a choice of normalization for $X$
 to prevent such spurious solutions as $(X,\tilde{A}) \rightarrow 0$
 and to prevent $X$ from becoming singular.

An error norm that doesn't rely on normalizing $X$ is $\|X^{-T} E X^{-1}\|_F$,
 with the corresponding minimization
\begin{equation} \min_{Y,\tilde{A}} {\left\| A - Y^T \tilde{A} Y \right\|_F}, \label{errmin2}\end{equation}
 where $Y \equiv X^{-1}$ is given the same local form as $X$.
This expression is less appealing because it is more nonlinear than Eq. (\ref{errmin})
 in that it contains $6^\textrm{th}$ order variable terms rather than just quartic terms.
However, it is a more direct minimization of the error perturbation
 that takes us from our approximate inverse to the exact inverse,
\begin{displaymath} A^{-1} = X \tilde{A}^{-1} X^T - (X \tilde{A}^{-1} X^T) (X^{-T} E X^{-1}) (X \tilde{A}^{-1} X^T) + O(\|X^{-T}EX^{-1}\|^2). \end{displaymath}
This error norm will be used for the remainder of this paper.

\subsection{Condition of the Transformation \label{illcond}}

Using a local transformation to remove matrix elements
 is only a specific application of a general ability
 to alter the values of matrix elements
 while preserving the sparsity pattern of a matrix.
We can consider a transformation to be part
 of a continuous family of transformations, $X(t)$ and $\tilde{A}(t)$,
 that begins at $t=0$ as the error free identity, $X(0)=I$ and $\tilde{A}(0)=A$,
 and ends at $t=1$.
We evolve from the error free transformation
 by following the minimum error transformations
 as we continuously turn on a non-negative constraint
 that enforces the final, restricted sparsity pattern at $t=1$,
\begin{equation} \min_{X(t),\tilde{A}(t)} \left( f_{error}[X(t),\tilde{A}(t)] + t \cdot f_{constraint}[X(t),\tilde{A}(t)] \right). \label{error_constrained}\end{equation}
Following this defined path of transformations,
 the constrained error norm in Eq. (\ref{error_constrained}) is non-decreasing with increasing $t$.
In order for the final transformation at $t=1$ to have a small error,
 the error must be small throughout the path
 and the Jacobian of error with respect to changes of $(X,\tilde{A})$
 must have an equally small near-null $(dX,d\tilde{A})$ component tangent to the path.
Correspondingly, we expect the condition number of the error minimization process
 to be inversely proportional to the minimum error attainable by the transformation.

\begin{figure} \begin{center}\includegraphics[width=200pt]{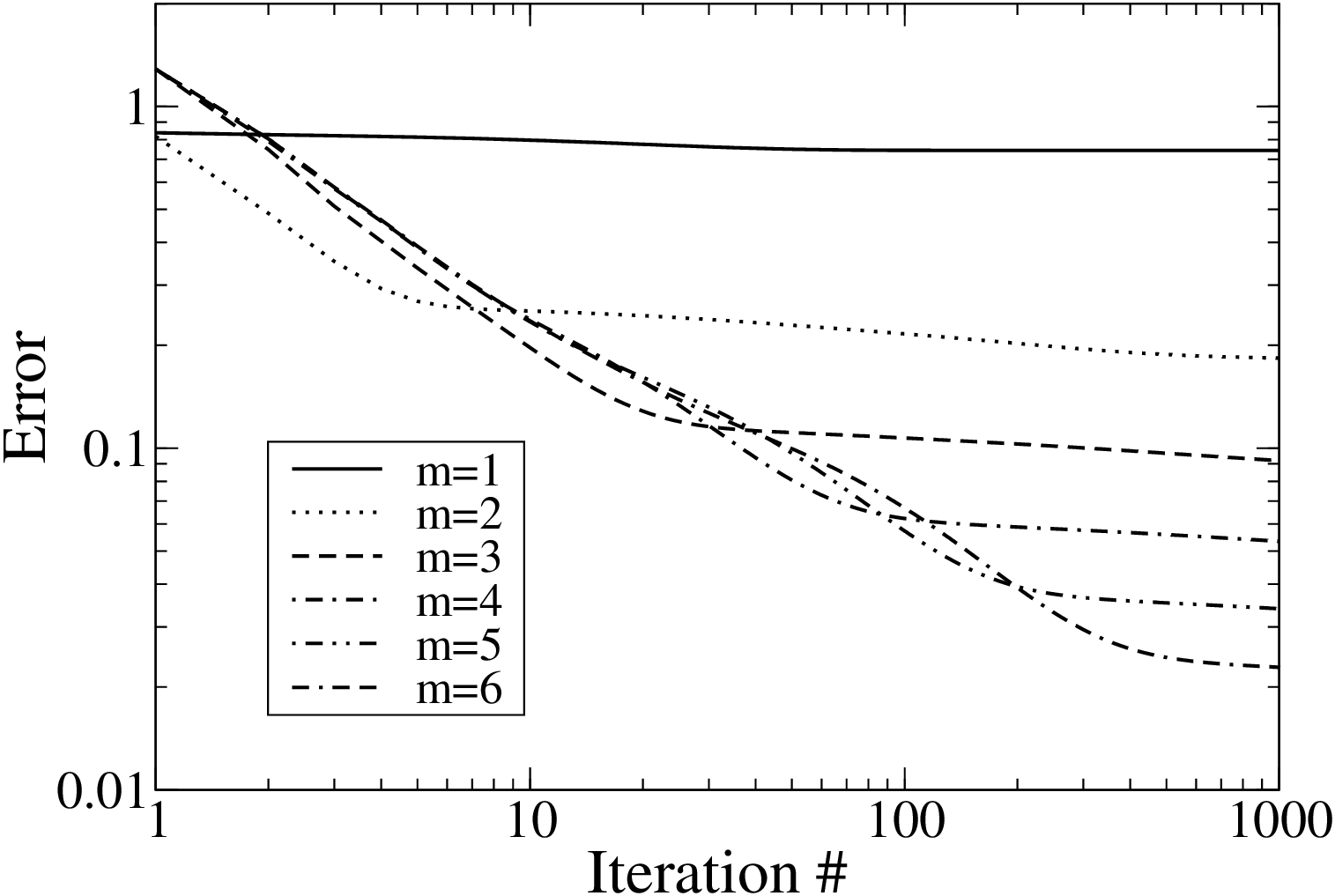}\end{center}
\caption{Convergence of error norm with steepest descent for $\lambda = 0$.} \label{fig_sd} \end{figure}

To illustrate the ill-conditioned nature of Eq. (\ref{errmin2}),
 we attempt to minimize it by following the negative gradient
 for a single transformation on our Helmholtz test problem at $\lambda = 0$.
We decouple one node on the interior of the grid
 without adding or subtracting any other terms from the sparsity pattern of $\tilde{A}$
 and the local region consists of all nodes within $m$ hops of the decoupled node.
We start from an initial guess of $Y=I$ and the nonzero terms of $\tilde{A}$ set to the corresponding values of $A$.
At each iteration, the gradient is calculated and the error norm is minimized in the direction of the gradient.
The error norm for the first 1000 iterations is plotted in Fig. \ref{fig_sd} for several values of $m$.
Only the $m=1$ case converges within 1000 iterations,
 but the expected trend of decreasing error and increasing condition number
 with increasing $m$ is readily apparent.
A tractable calculation of $Y$ and $\tilde{A}$ requires a more careful treatment of the ill-conditioned Jacobian.

\section{Linearized Approach to Local Coarsening\label{linmethod}}

Calculating and inverting the exact Jacobian of Eq. (\ref{errmin2})
 is impractical due to its size, ill-conditioning, and large null space.
The symmetric error matrix, $Y^T E Y$, contains $\frac{1}{2} n_L (n_L+1)$ elements to be minimized
 and the $Y$ and $\tilde{A}$ variables contain $(n_{\tilde{A}} + n_L n_I)$ independent unknowns.
For some local partitions, such as $m > 1$ in our test problem,
 there are more unknowns than matrix elements to be minimized,
 but the minimization is not underdetermined due to a large null space.
To alleviate these difficulties we separate the minimizations of $Y$ and $\tilde{A}$
 in an approximate way that leaves us with a well-conditioned problem in $Y$
 whose null space can be analytically removed
 and a much smaller, ill-conditioned problem in $\tilde{A}$.

We approximately linearize Eq. (\ref{errmin2}) by expanding $Y$ and $\tilde{A}$ in small changes,
 $Y \rightarrow (Y + dY)$ and $\tilde{A} \rightarrow (\tilde{A} + d\tilde{A})$,
 and keeping only terms to first order in $dY$ and $d\tilde{A}$ within the norm,
\begin{displaymath} \left\| Y^T (E-d\tilde{A}) Y - (P_I \tilde{A} Y)^T dY - dY^T (P_I \tilde{A} Y) \right\|_F. \end{displaymath}
Because of the restricted form of $Y$, $dY$ is an $n_I \times n_L$ matrix
 and $P_I$ is an $n_I \times n_L$ submatrix of identity.
This is not the correct way to linearize Eq. (\ref{errmin2}) -
there are additional linear terms proportional to $\|E\|$
 whose neglect leads to a linear convergence of the minimization -
 but this approximation leads to a greatly simplified solution.
The Frobenius norm is invariant with respect to orthogonal rotations of its operand,
 and we choose a particularly useful rotation
 consisting of the null space $\tilde{Q}$ and the spanned space $Q$ of $(P_I \tilde{A} Y)$.
Due to the $P_I$, the spanned space usually contains $n_I$ vectors,
 but it can contain less if $\tilde{A}$ is rank deficient.
If we apply the rotation to the error norm, we can write the norm squared as
\begin{equation} \begin{array}{l} \left\| Q^T Y^T (E - d\tilde{A}) Y Q - (P_I \tilde{A} Y Q)^T (dY Q) - (dY Q)^T (P_I \tilde{A} Y Q) \right\|_F^2 \\ + 2\left\| Q^T Y^T (E - d\tilde{A}) Y \tilde{Q} - (P_I \tilde{A} Y Q)^T (dY \tilde{Q}) \right\|_F^2 \\ + \hspace{390000 sp} \left\| \tilde{Q}^T Y^T (E - d\tilde{A}) Y \tilde{Q} \right\|_F^2. \end{array} \label{roterr}\end{equation}
The first two terms of Eq. (\ref{roterr}) can be canceled with a proper choice of $dY$,
\begin{displaymath} dY = \frac{1}{2}(P_I \tilde{A} Y)^{-T} \left( Y^T (E - d\tilde{A}) Y \right)(I + \tilde{Q} \tilde{Q}^T), \end{displaymath}
 where the inverse is a pseudo-inverse.
This leaves the third term to be minimized by $d\tilde{A}$,
\begin{equation} \min_{d\tilde{A}} {\left\| \tilde{Q}^T Y^T(E -  d\tilde{A}) Y \tilde{Q} \right\|_F}, \label{errmin3}\end{equation}
 which is overdetermined and ill-conditioned.

The approximately linearized Eq. (\ref{roterr}) has a significant null space,
 corresponding to additions to $dY$ of the form $(P_I \tilde{A} Y)^{-T}S(Q Q^T)$
 for any antisymmetric matrix $S$.
This null space has a size of $\frac{1}{2} n_I (n_I-1)$ for full rank $\tilde{A}$,
 which is large enough to account for Eq. (\ref{errmin2}) being overdetermined.

Solving Eq. (\ref{errmin3}) is the most difficult and expensive step of the error minimization.
The cost of an unstructured QR factorization of the problem is $O(n_B^2 n_{\tilde{A}}^2)$.
However, the system's matrix has some structure,
 it is a sum of two submatrices of the Kronecker product $(Y \tilde{Q})^T \otimes (Y \tilde{Q})^T$.
There are no existing structured QR factorization algorithms for this kind of matrix,
 but the structure allows for an efficient construction of the normal equations,
 which is a sum of two $n_{\tilde{A}} \times n_{\tilde{A}}$ submatrices of $(Y\tilde{Q}\tilde{Q}^T Y^T) \otimes (Y\tilde{Q}\tilde{Q}^T Y^T)$.
The cost of constructing and solving the normal equations is $O(n_{\tilde{A}}^3)$,
 which is an improvement over unstructured QR if $n_{\tilde{A}} \ll n_B^2$.
In our 2D example $n_{\tilde{A}} \sim n_B^2$
 making the order of complexity equal in both methods,
 but the normal equations are still faster due to a smaller prefactor.
The disadvantage of using the normal equations is the squaring of the condition number,
 which has a noticeable effect in the ill-conditioned, small error limit.

\subsection{Numerical Results\label{numerical}}

We return to the test problem at $\lambda = 0$,
 now using the linearized solution approach rather than following the gradient.
The same local region, $\tilde{A}$ sparsity pattern, and initial $Y$ and $\tilde{A}$ are used as in Section \ref{illcond}.
The $\tilde{A}$ minimization is performed using the normal equations
 which are solved using singular value decomposition (SVD) for testing purposes.
After each iteration the solution is updated,
 $Y \rightarrow (Y + \alpha dY)$ and $\tilde{A} \rightarrow (\tilde{A} + \alpha d\tilde{A})$,
 with $\alpha$ chosen to minimize the error norm.

\begin{figure} \begin{center}\includegraphics[width=300pt]{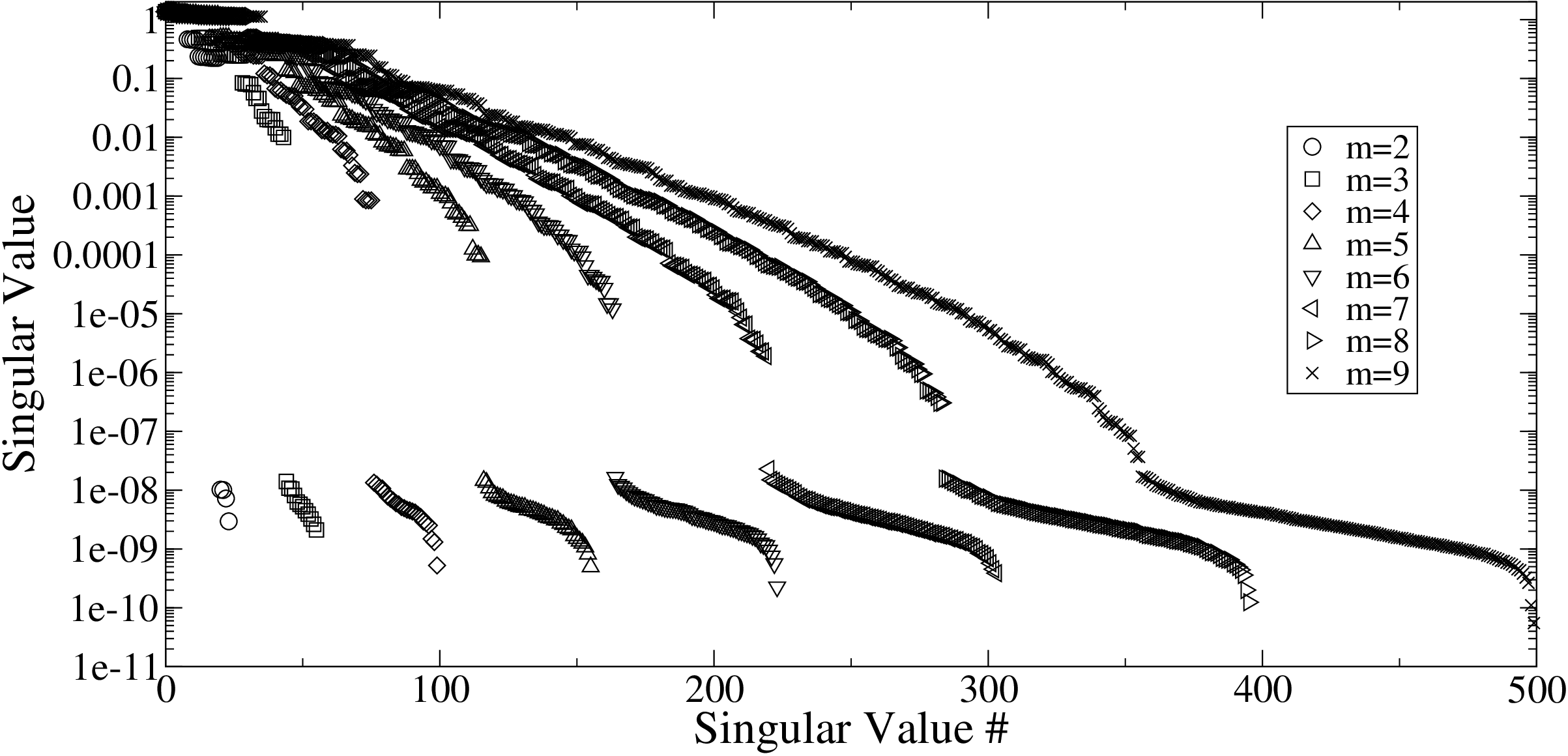}\end{center}
\caption{Singular values of Eq. (\ref{errmin3}) for $\lambda = 0$ (calculated from the normal equations).}
\label{fig_svd} \end{figure}

The SVD of Eq. (\ref{errmin3}), which is performed numerically on its normal equations,
 is shown in Fig. \ref{fig_svd} for the initial $Y$ and $\tilde{A}$.
An interesting feature of each SVD spectra is the null space of size $n_I$,
 resolved in this calculation to single precision, $10^{-8}$,
 relative to the largest singular value.
The null space corresponds to the set of local transformations that exactly preserve sparsity
 and in this case diagonal scaling of the interior block, $(Y,\tilde{A}) \rightarrow (D^{-1}Y,D\tilde{A}D)$.
A change in diagonal scaling doesn't effect the chosen error norm from Eq. (\ref{errmin2})
 and correspondingly the error term in Eq. (\ref{errmin3})
 is orthogonal to the null space within machine precision.
The contribution to $d\tilde{A}$ from the null space should be zero,
 and since it can be clearly distinguished from the gap in the spectrum,
 we can simply ignore the null space component.
The smallest singular value of the rest of the spectrum
 shows an exponential decay with respect to $m$,
 which suggests an exponential decay of the minimum error according to the argument in Section \ref{illcond}.
The limiting effects of finite precision are clearly visible
 in the vanishing of the gap between the null and spanned space for $m \ge 9$.

\begin{figure} \begin{center}\includegraphics[width=200pt]{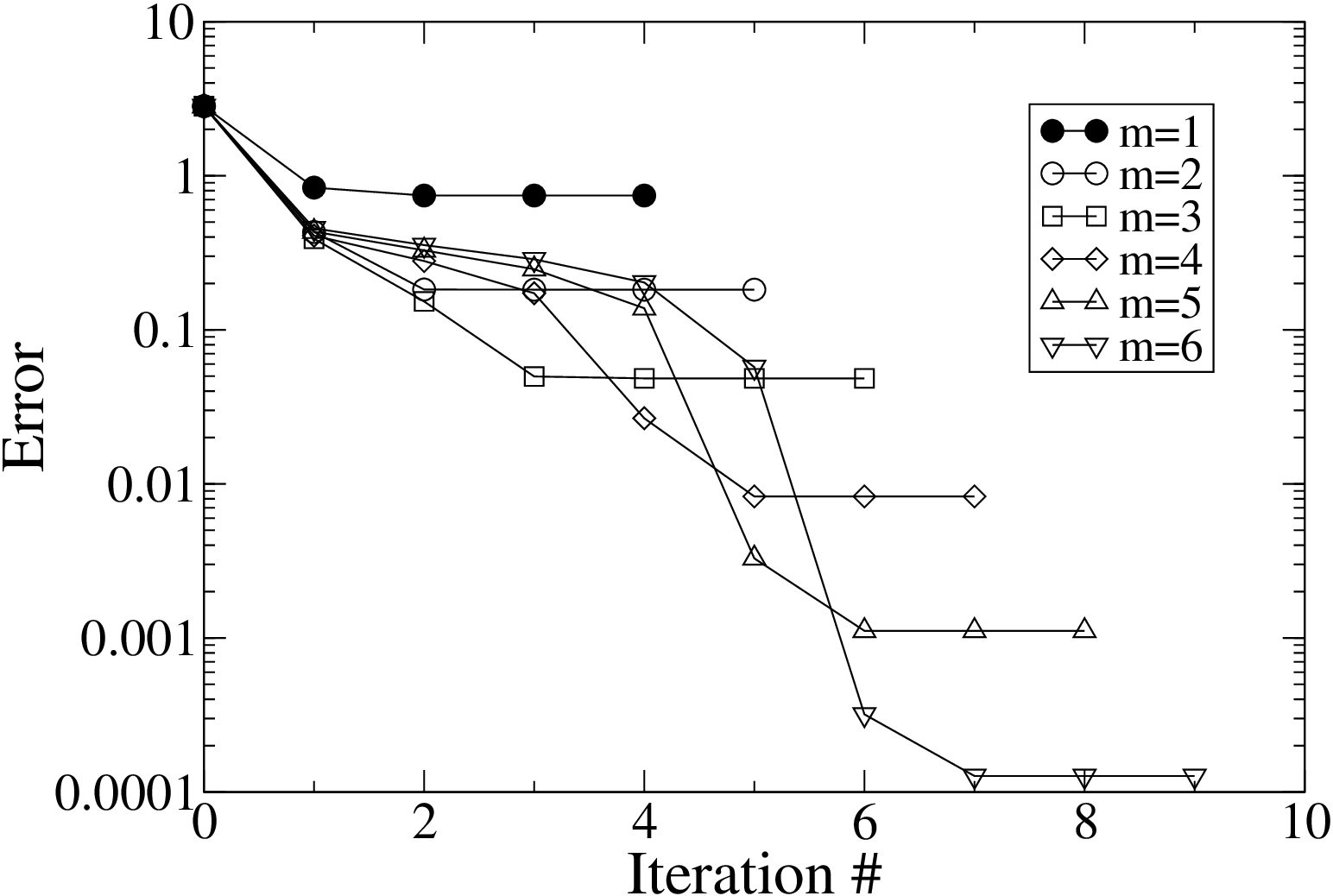}\end{center}
\caption{Convergence of error norm with linearized solutions for $\lambda = 0$.} \label{fig_conv} \end{figure}

\begin{figure} \begin{center}\includegraphics[width=200pt]{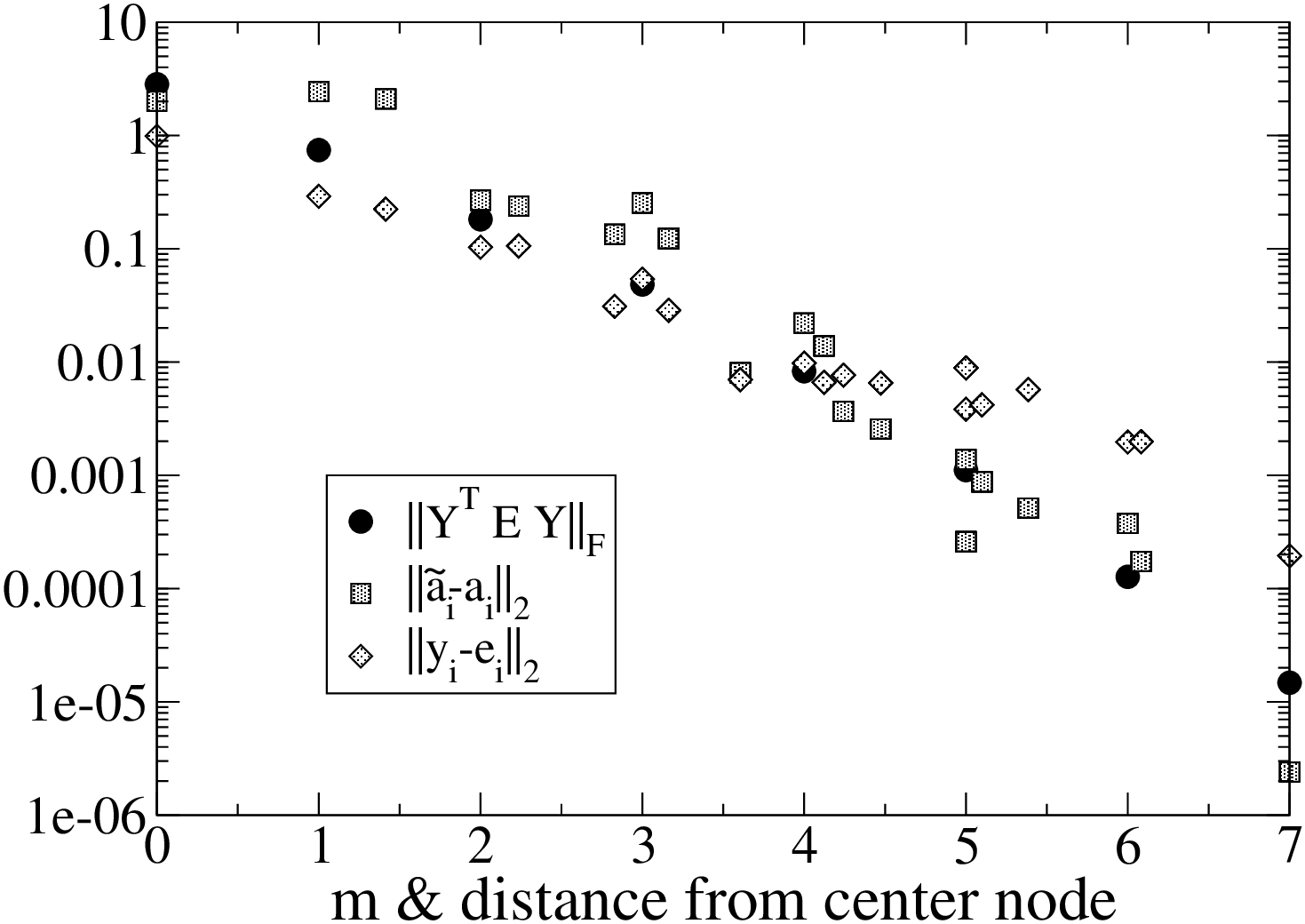}\end{center}
\caption{Spatial decay comparison, minimum error norm for different $m$ values and column norms from the $m=7$
minimum error solution.} \label{fig_space} \end{figure}

The convergence of the linearized solution approach is shown in Fig. \ref{fig_conv}.
Each calculation takes approximately $m$ steps to converge,
 which signifies the success of our approximate inverse Jacobian
 in capturing the ill-conditioned aspects of the problem.
The error exponentially decays with $m$ as expected from the spectrum of Eq. (\ref{errmin3}).
This spatial decay of error can be related
 to a spatial decay of $\tilde{A}$ to $A$ and $Y$ to $I$
 by taking the error to be caused by the truncation of some dense exact $Y$
 to $I$ outside a local region.
The relation of the decays can be seen in Fig. \ref{fig_space},
 where the error norm as a function of $m$ is plotted against
 the deviations of $Y$ and $\tilde{A}$ from $I$ and $A$ measured by column
 and plotted by the geometric distance on the 2D grid of the associated node from the central, decoupled node.
Since $Y$ and $\tilde{A}$ are only defined up to a diagonal scaling of the interior block,
 the rows of $Y$ are normalized to a 2-norm of one to make them unique.

\begin{figure} \begin{center}\includegraphics[width=200pt]{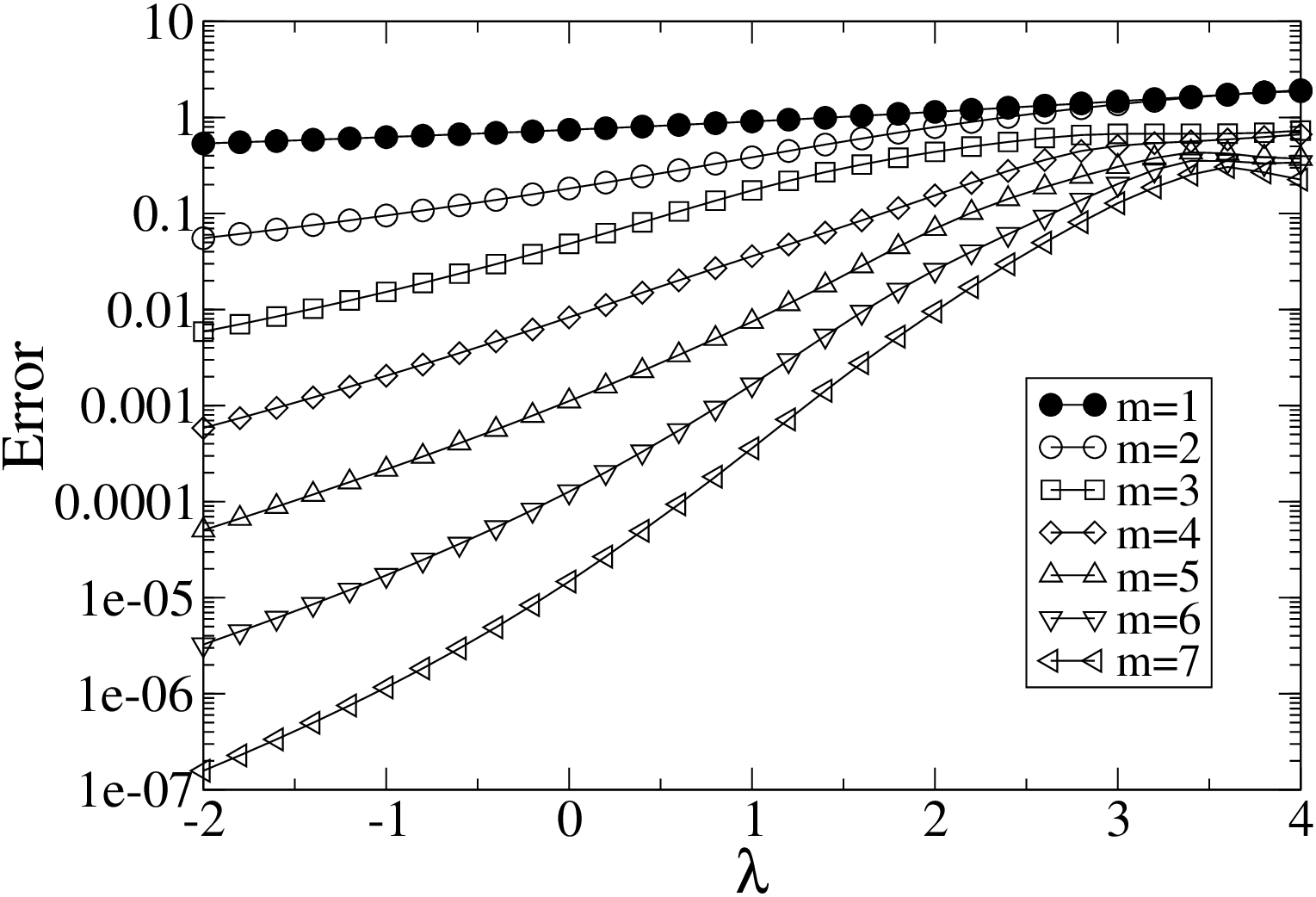}\end{center}
\caption{Converged error norm versus $\lambda$ using linearized solutions.} \label{fig_lambda} \end{figure}

\begin{figure} \begin{center}\includegraphics[width=200pt]{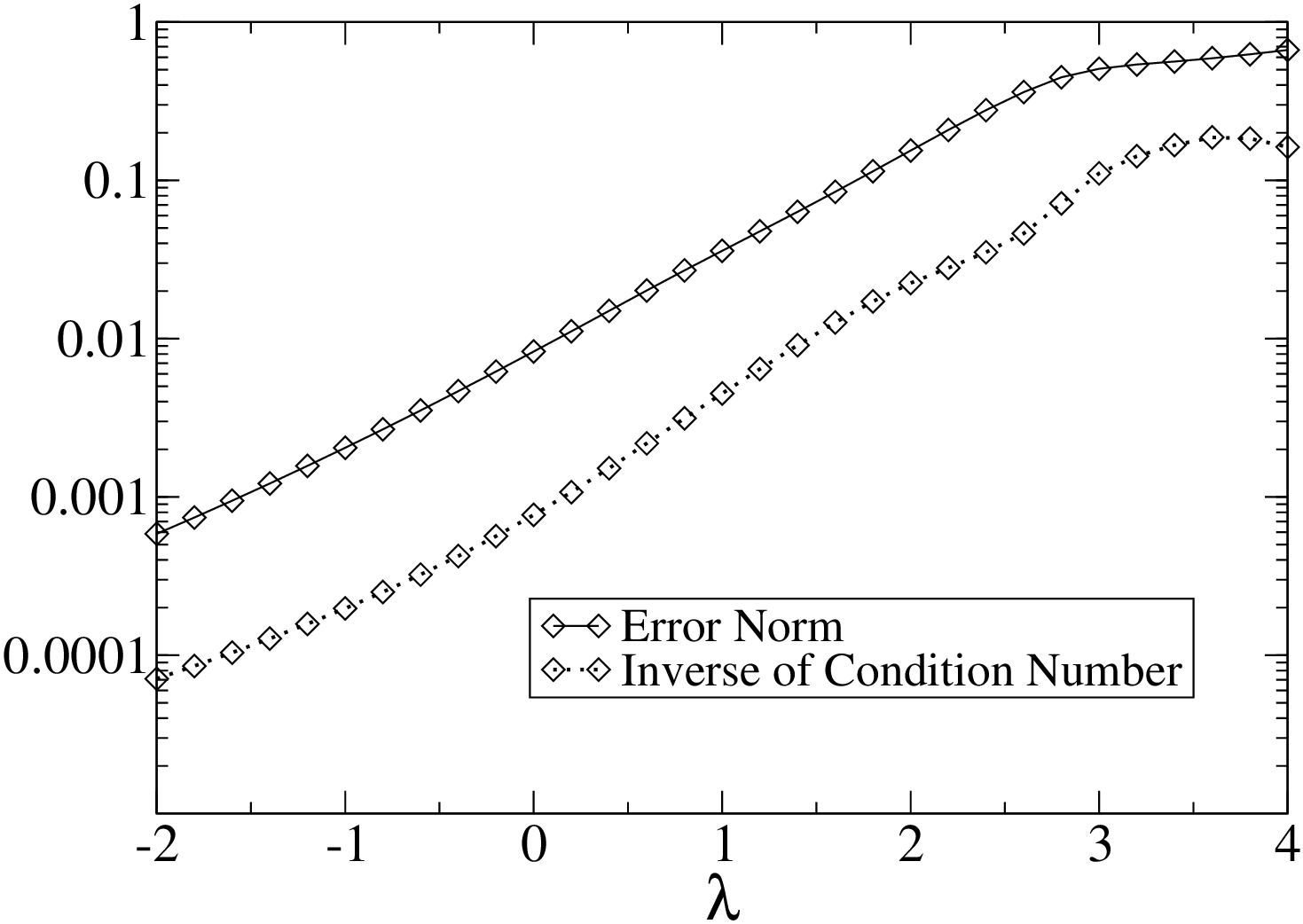}\end{center}
\caption{The error norm and condition number of Eq. (\ref{errmin3}) at convergence for $m=4$.}
\label{fig_errcond} \end{figure}

We next try the method on the more interesting $0 < \lambda < 8$ case,
 though only the $0 < \lambda \le 4$ range needs to be tested
 as the matrix for $\lambda = 4 - x$ can be mapped to $\lambda = 4 + x$ with a diagonal scaling.
The converged error norm for multiple values of $\lambda$ and $m$ are plotted in Fig. \ref{fig_lambda}.
The condition number of the row normalized $Y$ is less than twelve for all calculations performed.
The decay of $A^{-1}$ off the diagonal is exponential in geometric distance for $\lambda < 0$,
 but this qualitative change in behavior from $0 \le \lambda \le 8$
 doesn't cause any kinks in the error at $\lambda = 0$.
We observe that the exponential decay rate of error with $m$ is approximately proportional to $|\lambda-4|$.
Near the $\lambda = 4$ point,
 the exponential error decay appears to break down
 leaving an error with an $m$-dependence proportional to the logarithmic decay of off-diagonal elements of $A^{-1}$.
The most obvious matrix property to attribute to the loss of decay near $\lambda = 4$
 is the vanishing of the diagonal elements of $A$.

The inverse proportionality between the condition number of the non-null subspace of Eq. (\ref{errmin3})
 and the minimum error norm continues to hold as a function of $\lambda$ as shown in Fig. \ref{fig_errcond}.
The condition number plotted is calculated at the converged $(Y,\tilde{A})$ value,
 but the condition number varies very little between iterations
 and it is within a factor of two of the condition number calculated from the initial $(Y,\tilde{A})$ guess.

The loss of exponential error decay as $\lambda \rightarrow 4$
 signifies the disappearance of locally removable degrees of freedom
 from a model restricted in form by the restriction on the sparsity pattern.
As $\lambda \rightarrow 4$ the wavelength of oscillations in the Helmholtz equation
 approaches four times the grid spacing, a high frequency limit where multigrid also fails.
For the multigrid approach to continue into this limit,
 the solution must be decomposed into a sum of envelope functions
 times oscillatory solutions with wavevectors in various directions \cite{waveAMG}.
This is a transformation from a scalar differential equation
 to a vector differential equation,
 which can't be represented by $A \rightarrow \tilde{A}$
 unless the sparsity pattern of $\tilde{A}$ is allowed to fill in somewhat.
Error decay is restored for $\lambda > 4$ only because
 the discretization of the Helmholtz equation breaks down
 and the correct high frequency oscillations are no longer present in the matrix problem.

\section{Choice of Sparsity Pattern \label{blocksection}}

For all the tests performed in Section \ref{linmethod}
 we strictly prevented fill-in in transforming from $A$ to $\tilde{A}$,
 but it is only really necessary to control fill-in
 enough to preserve the scalability of the factorization.
It can be beneficial to add nonzero matrix elements to $\tilde{A}$
 because that increases the number of degrees of freedom in the error minimization, Eq. (\ref{errmin2}),
 and can reduce the minimum error norm.

An important reason for filling in $\tilde{A}$
 is to prevent the removal of a node from
 breaking the global connectivity of a problem.
The simplest case of this is a tridiagonal matrix,
 which can be associated with a problem on a 1D grid.
If a node is removed from the grid without filling in the matrix,
 then the grid will be split in half.
The associated transformation would have to contain
 all the response of each half of the grid on the other
 and cannot in general be accurately made local.

One simple way to avoid changing the connectivity of a problem
 is to aggregate nodes together into supernodes
 where all the member nodes share all the connections of other member nodes.
Once a supernode is formed, the decoupling of a node in the supernode
 won't break any connectivity as long as one node remains within the supernode.
The supernode concept has been used before in sparse Gaussian elimination for efficiency reasons \cite{SUPERLU}
 to allow for the use of dense matrix operations in inner loops,
 but here it serves a more fundamental purpose.
The larger the supernodes are made, the more filled in the matrix will be,
 and the error norm will have a decreasing minimum with fixed local region size.
If the supernodes are made large enough,
 then Gaussian elimination steps can be performed without additional matrix filling
 before the more expensive algebraic coarsening procedure
 in a possibly more efficient hybrid approach.
For our example on a 2D grid, the grid of nodes can be made a grid of supernodes,
 which can be interpreted as a discretization of a vector differential equation
 where the number of vector components is the size of the supernode.

\begin{figure} \begin{center}\includegraphics[width=200pt]{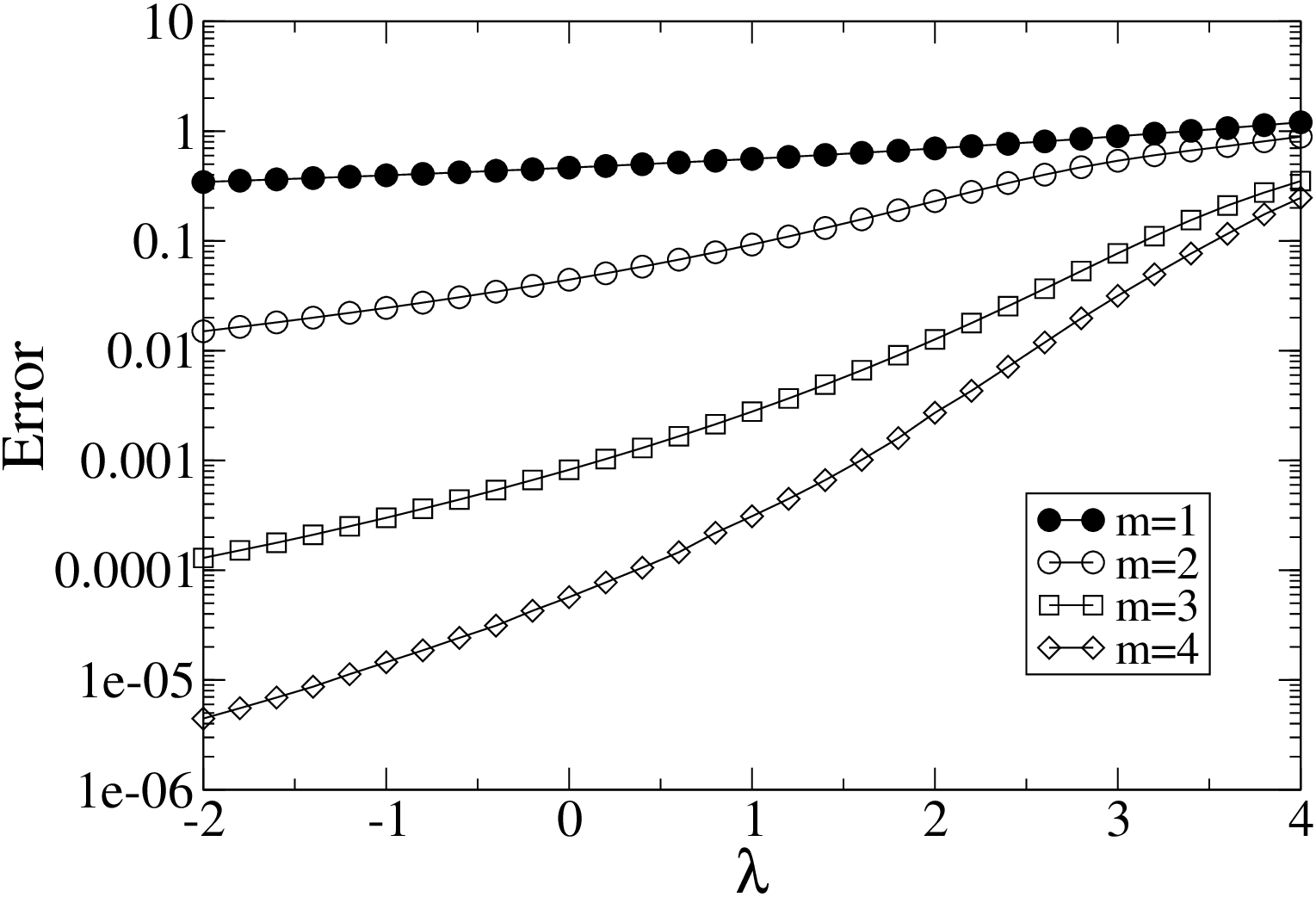}\end{center}
\caption{Converged error norm for a grid of $p=2$,$q=1$ supernodes.} \label{fig_block} \end{figure}

\begin{figure} \begin{center}\includegraphics[width=200pt]{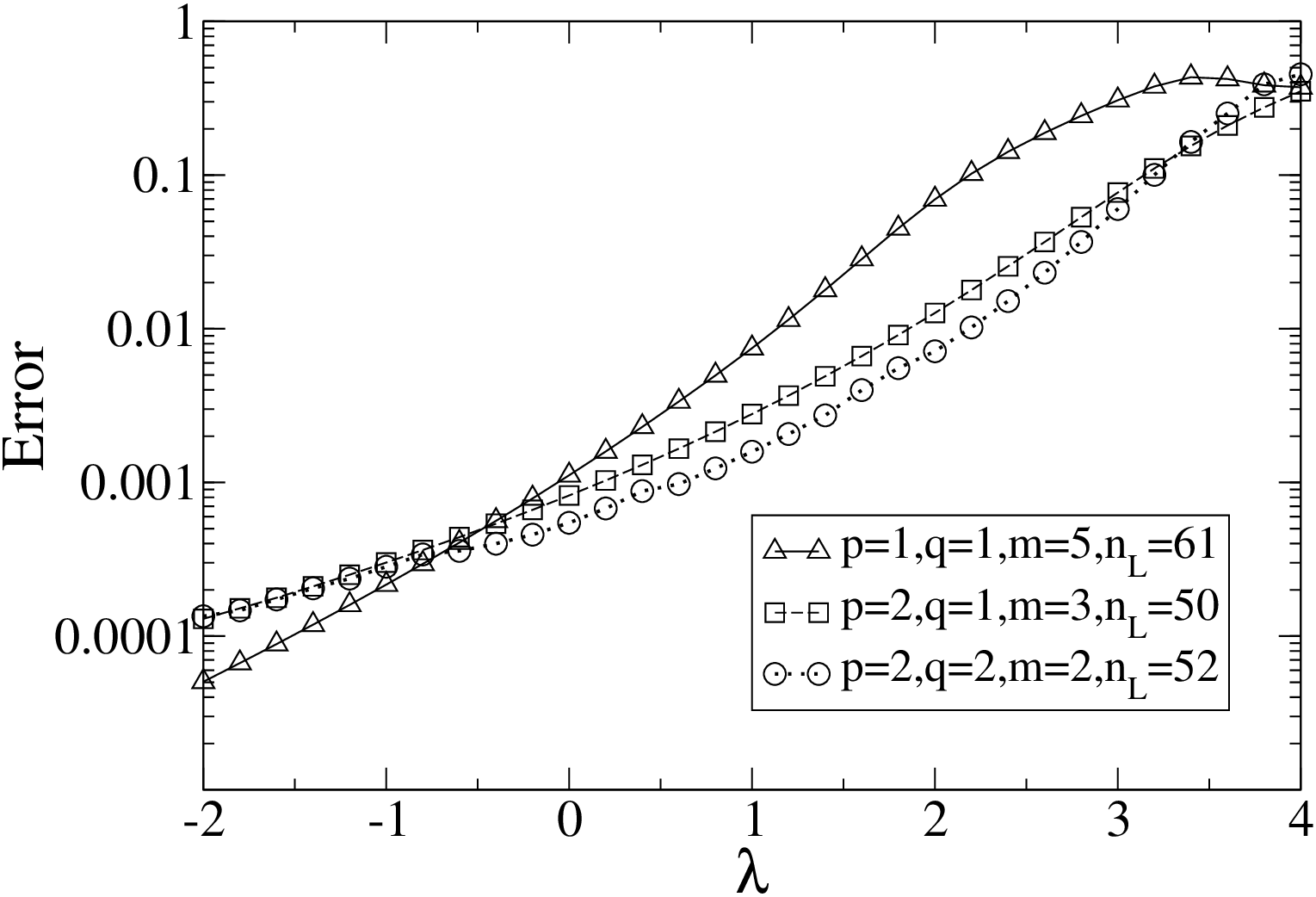}\end{center}
\caption{Converged error norm for different supernodes and similar $n_L$.} \label{fig_multiblock} \end{figure}

We again return to the test problem, now with a 2D grid of supernodes
 constructed by merging $p \times q$ rectangles of neighboring nodes.
A local transformation is performed to remove one node from one supernode
 and the local region is chosen to include all supernodes within $m$ hops
 of the removed node.
The converged error norm for $p=2$, $q=1$ is plotted in Fig. \ref{fig_block}
 and the important difference with Fig. \ref{fig_lambda}
 is the error seems to continue to decay exponentially in $m$ near $\lambda = 4$
 rather than stagnate at $\lambda \cong 3.5$.
A comparison between three different supernode sizes is plotted in Fig. \ref{fig_multiblock}.
For similar $n_L$ all errors are roughly the same in the non-oscillatory regime, $\lambda \le 0$,
 and at the maximally oscillatory point, $\lambda = 4$,
 while the larger supernodes' errors are smaller in the intermediate oscillatory regime, $0 < \lambda < 4$.
This result suggests that supernodes are useful
 for increasing the rate of error decay,
 but if the $\lambda = 4$ case can indeed be improved,
 larger supernodes are required.

The arguments made in Ref. \cite{waveAMG} suggest that at least an eight wave expansion
 is required for an efficient solution of the $\lambda = 4$ case,
 which might be properly captured by $p=2,q=4$ or $p=3,q=3$.
However, both cost and conditioning are a barrier
 to the current approach to calculating these transformations.
The cost of solving the normal equation version of Eq. (\ref{errmin3})
 for fixed local region size scales as $O(p^3 q^3)$.
The conditioning of Eq. (\ref{errmin3}) remains
 inversely proportional to the minimum error norm,
 but the constant of proportionality is observed
 to change substantially with $p$ and $q$,
 causing calculations to be more ill-conditioned with the same minimum error norms.

\section{Conclusions}

We have studied the possibility of factoring sparse matrices
 by means of local sparsity-preserving transformations
 with numerical tests of a single transformation.
Intermediate stages of such a factorization require
 transformations to be performed on matrices of similar sparsity
 but with different values of their matrix elements,
 which was examined here in a simple, artificial manner
 by varying the frequency of the Helmholtz equation.
Qualitatively, we expect the intermediate, ``coarsened"
 matrices to still represent the Helmholtz equation
 with the frequency scaled to represent a change of length scale.
The absence of local degrees of freedom
 for the $\lambda = 4$ case in Section \ref{numerical}
 suggests that this interpretation fails when the wavelength becomes proportional
 to the grid spacing.
To continue to remove local degrees of freedom beyond this frequency,
 it becomes necessary to allow the coarsened matrices to take a more general form.

The sparse linear solver methodology presented in this paper
 has demonstrated a behavior distinct from both direct and iterative solvers.
The success of direct solvers is dependent on the filling in of the matrix
 in the intermediate stages of factorization, which is a graph theoretic property
 and is controlled by the order in which nodes are factored.
The success of iterative solvers is dependent on the condition number of the matrix
 and is controlled by preconditioning a problem to reduce the condition number.
Here the determining characteristic of how costly it is to solve a matrix
 is the decay of error of a transformation with respect to local region size
 and can be controlled by changing the local region or sparsity pattern.
Matrix fill is no longer a problem as it is strictly controlled,
 and the error decay is a local property completely independent of the global spectrum
 and conditioning of the matrix.

There remain technical difficulties with calculating local sparsity-preserving transformations
 that must be resolved before a practical linear solver can be implemented with them.
The most important problem is determining whether or not a
 well-conditioned process exists for calculating local transformations.
The ill-conditioning is associated with minimizing an error norm,
 and a method based on additional criteria might precondition the process.
Another important problem is understanding what properties of a matrix
 and sparsity pattern determine the rate of decay of error.
This is needed to determine precisely when sparsity patterns should be changed
 and how they should be changed to make calculations most efficient.
Once the issues associated with single transformations are resolved,
 there is the additional problem of choosing the ordering of transformations.
The ordering can determine error decay rates of successive transformations,
 error propagation during factorization,
 and the amount to which the process can be parallelized.
The simple answer at least for the purpose of parallelization
 is to choose as many transformations as possible on disjoint local regions
 to maximize the number of concurrent calculations of local transformations.

\section*{Acknowledgements}

This work was supported by National Science Foundation Grant No.
DMR04-39768 and by the Director, Office of Science, Office of Basic
Energy Sciences, Division of Materials Sciences and Engineering, U.S.
Department of Energy under Contract No. DE-AC03-76SF00098.

I would like to thank Professor Marvin Cohen for guidance and support.

I also thank Jay Deep Sau, David Bindel, Oren Livne, and Raim Tavisseur
 for numerous useful discussions.


\end{document}